%% file: elaz-remark.tex
\documentclass[12pt]{amsart}
\usepackage{amsmath}
\usepackage{amssymb}
\usepackage{amsfonts}
\usepackage{amsthm}
\usepackage{verbatim}
\usepackage{amscd}
\usepackage{cite}
\usepackage{leftidx}
\usepackage{enumerate}
\usepackage{txfonts}
\usepackage{manfnt}
\usepackage{amscd}
\usepackage[mathscr]{eucal}
\usepackage{hyperref}
\usepackage{graphicx}
\usepackage{palatino}

\input{myheader3.tex}

\newtheorem{thm}{Theorem} % [section]

\newtheorem*{thm*}{Theorem}
\newtheorem*{prop*}{Proposition}
\newtheorem{cor}[thm]{Corollary}
\newtheorem*{cor*}{Corollary}

\newtheorem*{claim*}{Claim}

\theoremstyle{remark}

\newtheorem{rem}[thm]{Remark}
\newtheorem{crit-rem}[thm]{Critical remark}

\newtheorem*{remarks*}{Remarks}

\newtheorem{example}[thm]{Example}
\newtheorem*{example*}{Example}

\newtheorem*{defn*}{Definition}

\begin{document} 
\title{Asymptotic syzygies grow
 exponentially:\\
A Remark on a paper of Ein-Lazarsfeld}
\author{Daniel Chun}
\author {Ziv Ran}
{}
%\author{}
%\Large

%\thanks{\raggedright{

%Partially supported by NSA Grant MDA904-02-1-0094} }
%\thanks{arxiv.org/1308.2442}
\date {\today}% \enddate

%\affil University of California, Riverside\endaffil

\address {\nl UC Math Dept. \nl
Big Springs Road Surge Facility
\nl
Riverside CA 92521 US\nl 
ziv.ran @  ucr.edu\nl
\url{http://math.ucr.edu/~ziv/}
}

%\email {ziv.ran @ucr.edu}
% \subjclass[2010]{14J40, 32G07, 32J27, 53D17}
%\keywords{}

\begin{abstract}
%We compare Koszul cohomology of a variety and a subvariety
%via certain filtrations on Koszul complexes.
For all $2\leq q\leq \dim(X)$ and most relevant $p$ values, the dimension
of the asymptotic
Koszul cohomology group $K_{p,q}(X, B;L_d)$ grows exponentially with $d$.
\end{abstract}
\maketitle

%\section{Introduction}
While interest in the syzygies of a projective variety $X$ goes back to Hilbert,
it was Mark  L. Green's seminal paper \cite{green-koszul} that popularized the notion
that the syzygies of $X$ can be viewed as a type of cohomology,
called Koszul cohomology, attached to the pair $(X, L)$ where $L$ is a line
bundle on $X$. Since then, chiefly under
the leadership of Robert Lazarsfeld and Claire Voisin, numerous results
have been established on the computation or nonvanishing
 of syzygies in many cases.
 Recently
particular attention has been focused on the case
of \emph{asymptotic} Koszul cohomology groups
$K_{p,q}(X, B;L_d)$ where $L_d$ is sufficiently ample
 \cite{el}. Some of the results are stated in \cite{ein-lazarsfeld-yang}. \par
Larzarsfeld and Ein \cite{el} have shown, for a large range of $p$ and $q\geq 2$ that these groups
don't vanish. Our purpose here is to show that the construction of \cite{el} can be pushed to
yield lower bounds on the dimension of these
groups which are exponential in $d$. Our argument differs from \cite{el} mainly
 in the elimination of the induction, i.e. 'drilling down to the base case'.\par
 After this was written, Lazarsfeld communicated that he and Ein had been aware
 in a general way of
 similar bounds, based on an argument different from the one given here, which was never
 written down.  Nonetheless, the asymptotic behavior of the
 dimension of the $K_{p,q}(X, B; L_d)$ is stated
 an an open problem in \cite{el}. See \cite{el-sav} for further related results and problems.

%Our purpose here is to work out a local cohomology version of (asymptotic) Koszul 
%cohomology, with support in a subvariety $Y\subset X$.
%This takes the form of some filtrations on Koszul complexes on $X$,
%whose graded pieces are related either to $(Y, L|_Y)$ or to the projection
%of $X$ from $Y$ (or equivalently, the blowup of $X$ in $Y$
%and the sections of $L$ vanishing on $Y$).
%As a consequence we will prove a result on nonvanishing of certain asymptotic
%Koszul cohomology groups first obtained with a different proof by Ein and
%Lazarsfeld \cite{ein-lazarsfeld-inventiones}. Our proof
%is inspired by the argument in \cite{ein-lazarsfeld-inventiones}
%which also involves a (different) comparison with a different,
%trickier type of
%subvariety. Here we will use a simple complete intersection.
\section*{Notations and conventions}
Let $L$ be a line bundle and $B$ a coherent sheaf on a variety $X$. 
A \emph{linear system} on $X$, belonging to $L$, is by definition a 
vector space $V$ endowed with
a map, not necessarily  injective, $V\to H^0(X, L)$.
This induces  a map of coherent sheaves $V\otimes\O_X\to L$. 
 In our applications,
we will initially take $V$ complete, i.e. $V=H^0(X, L)$, but various
associated linear systems, e.g. restrictions on subvarieties, need not be complete.

Given these data, the associated
\emph{Koszul complexes} are by definition the complexes
\eqspl{koszul}{
C_{p, c}(X, B, L, V):\ \  \bigwedge\limits^pV\otimes(B+(c+1)L)\to\bigwedge\limits^{p-1}V\otimes
(B+(c+2)L)\to...
}
where we use customary additive notation: $B+aL:=B\otimes L^{\otimes a}$.
We will often suppress $X, L$ when understood and denote this by $C_{p, c}(B, V)$
or $C_p(B,V)$ when $c=0$.
When $V$ is base-point free, the
complex \eqref{koszul} admits local homotopy operators, hence it is
locally acyclic, i.e. exact in strictly positive degrees.
We will denote the kernel of \eqref{koszul}, i.e. $\H^0(C_p(B,V))$, by $M_p(B, V)$.
In our applications, $L$ will have the form
\[L=L_d:=dA+P\]
with $A$ ample and $d$ as large as desired. Then the members of the
Koszul complex are themselves acyclic sheaves thanks to Serre vanishing.
Whenever the Koszul complex is locally 
acyclic and consists of acyclic sheaves, 
we have
\eqspl{}{
\HH^i(C_{p,c}(B, V))=H^{i}(M_{p,c}(B, V))=K_{p-i,c+1+i}(X, L, B, V), i>0.
}
It will be convenient to use the following
terminology: a nonnegative complex $C^\bullet$ of sheaves
is \emph{globally
acyclic} if it is locally acyclic and consists of acyclic sheaves. In this
case, the hypercohomology can be identified
with either the cohomology of the $\H^0$ sheaf, or with the
cohomology of the associated $H^0$ complex. Thus, we have
\eqspl{}{
\HH^i(C^\bullet)= H^i(\H^0(C^\bullet))=H^i(H^0(C^\bullet)).
}
If $C^\bullet$ is globally acyclic
and these groups are moreover trivial for all $i>0$, $C^\bullet $ is said
to be \emph{hyper-acyclic}.\par
We will often use $H^\bullet$ rather than $\HH^\bullet$ to denote hypercohomology of
a complex of sheaves, when there is no confusion.

\section*{Exponenetial Ein-Lazarsfeld}
Let $Y\subset X$ be a smooth subvariety of dimension $m=n-c$ and consider a subspace $V_0\subset
V^+_0:=V\cap H^0(\I_Y(L))$ whose zero-locus is $m$-dimensional.
In our application, we will take $V=H^0(X, L)$.
Let $V_Y=V/V_0$, which admits a natural map to $L_Y:=L.\O_Y$. 
Suppose $p\geq v_Y:=\dim(V_Y)$.
Then the complex $C_{v_Y}(Y, B_Y, L_Y, V_Y)$ is exact,
hence so is $C_{v_Y}(Y, B_Y, L_Y, V_Y)\otimes \wedge^{p-v_Y}V_0$.
Now there is a
a natural quotient map 
\[\wedge^iV\to\wedge^{i-v_Y}V_0\otimes
\wedge^{v_Y}V_Y\]
(where the target is taken to be $\wedge^{v_Y}V_Y$ if $i=v_Y$ and zero if
$i<v_Y$). This induces a
 natural termwise surjective
map  
\[C_p(X, B, L, V)\to C_{v_Y}(Y, B_Y, L_Y, V_Y)\otimes \wedge^{p-v_Y}V_0 .\] 
The kernel of the latter map, denoted
$C^0_p(X, B, L, V)$, is hence a subcomplex 
quasi-isomorphic to $C_p(X, B, L, V)$. We call it the \emph{adapted} Koszul complex
(with respect to $Y$). Note that there is a natural map
\eqspl{x}{C^0_p(X, B, L, V)\to\wedge^{v_Y}V_Y\otimes\wedge^{p-v_Y}V_0\otimes\I_{Y/X}(B+L).
}
In the sequel we will use for $Y$ the subvariety used by Ein-Lazarsfeld \cite{el}, 
which in turn
is based on a construction of Eisenbud et al. \cite{eisenbud-restricting}. 
For a very ample divisor
$H$ fixed independently of $d$,
this $Y$
is of type 
\eqsp{H_1\cap...\cap H_{c-2}\cap D_1\cap D_2,
H_i\in |H|, D_1\in|H+B-K_X|, D_2\in|L-(c-1)H|
.}

$Y$ has dimension $m=q-2$ and
satisfies 
\eqspl{vanishing}{
h^{q-1}(\I_{Y/X}(B+L))=1; h^{q-1+i}(\I_{Y/X}(B+L+kH))=0, k>i\geq 0. 
}
Now let $\bar X$ be a smooth subvariety of $X$ transverse to (or disjoint from) $Y$ and let $\bar L, \bar B$ etc
denote restrictions on $\bar X$. Then for any divisor $B'$  we have an analogous map
\eqspl{xbar}{C^0_p(\bar X, \bar B', \bar L, V)\to\wedge^{v_Y}V_Y
\otimes\wedge^{p-v_Y}V_0\otimes\I_{\bar Y/\bar X}(\bar B'+\bar L).
}
Let $V_K, \bar V$ denote respectively the kernel and image of
the restriction map $V\to H^0(\bar L)$.
Then we have an exact diagram
\eqspl{v}{
\begin{matrix}
&0&&0&&0&\\
&\downarrow&&\downarrow&&\downarrow&\\
0\to&V_0\cap V_K&\to&V_K&\to&V_{K,Y}&\to 0\\
&\downarrow&&\downarrow&&\downarrow&\\
0\to&V_0&\to&V&\to&V_Y&\to 0\\
&\downarrow&&\downarrow&&\downarrow&\\
0\to&\bar V_0&\to&\bar V&\to&\bar V_Y&\to 0\\
&\downarrow&&\downarrow&&\downarrow&\\
&0&&0&&0&.
\end{matrix}
} 
\par
Now take $\bar X$ general of the form
\[\bar X=H'_1\cap...\cap H'_{m+1}, H'_i\in|H|\]
where $H$ is a fixed very ample divisor as above. Note that
a Koszul resolution of $\O_{\bar X}$ on $X$ remains exact
upon tensoring with $\I_Y$ and its locally free twists, and yields a torsion-free (not locally free)
Koszul resolution of $\O_{\bar X}$, since $\bar Y=Y\cap\bar X=\emptyset$ by
transversality. 
Thus we have exact
\small
\eqspl{koszul-torsionfree}{
\I_{Y/X}(B+L)\to(m+1)\I_{Y/X}(B+L+H)\to...\to\I_{Y/X}(B+L+(m+1)H)
\to\O_{\bar X}(B+L+(m+1)H).
}
\normalsize
Then setting $B'=B+(m+1)H$, consider an analogous 
torsion-free Koszul resolution of
$\wedge^{p-v_Y}V_Y\otimes\I_{\bar Y/\bar X}(\bar B'+\bar L)=
\wedge^{p-v_Y}V_Y\otimes \O_{\bar X}(\bar B+\bar L)$ , 
together with a compatible one for 
$C^0_p(\bar X, \bar B', \bar L, V)=C_p(\bar X, \bar B', \bar L, V)$.
This takes the form of a morphism of complexes
\tiny
\eqspl{map}{
\begin{matrix}
C^0_p(X, B, L, V)&\to...\to&(m+1)C^0_p(X, B+L+(m+1)H)&\to&C_p(\bar X, \bar B', \bar L, V)\\
\downarrow&&\downarrow&&\downarrow\\
\wedge^{v_Y}V_Y\otimes\wedge^{p-v_Y}V_0\I_{Y/X}(B+L)&\to...\to& (m+1)\wedge^{v_Y}V_Y\wedge^{p-v_Y}V_0\otimes\I_{Y/X}
(B+L+(m+1)H)&\to&\wedge^{v_Y}V_Y\wedge^{p-v_Y}V_0\otimes \O_{\bar X}(\bar B'+\bar L)
\end{matrix}
}
\normalsize
Pick a natural number $\bar p\leq v_Y$ independent of $d$, e.g. $\bar p=1$,
 and choose $V_0 $
so that $\dim(\bar V_Y)=\bar p$. 
For example, we can choose $V_0$ generated by $n$ general elements of 
$V_0^+$ (which cut out $Y$), plus the inverse image in $V_0^+$
of a general codimension-$(\bar p+n)$ subspace of the image of
$V_0^+\to\bar V$ (cf. \eqref{v}).
Note that in that case
$\dim(\bar V_0)$ is asymptotically $d^m$ for large d.
%hence $\dim(V_0\cap V_K)$ is asymptotically $d^n$. 
Then
the rightmost column of \eqref{map}
fits in a diagram, in which the upper horizontal arrows represent
a direct summand inclusion (with a suitable shift):
\eqspl{xbar1}{
\begin{matrix}
C_{\bar p}(\bar X, \bar B', \bar L, \bar V)\otimes\wedge^{p-\bar p}V_K&
\leftrightarrows& 
C_p(\bar X, \bar B', \bar L, V)\\
\downarrow&&\downarrow\\
\wedge^{v_Y}V_Y\otimes \wedge^{p-v_Y}(V_K\cap V_0)\otimes\O_{\bar X}
(\bar B'+\bar L)&\hookrightarrow&\wedge^{v_Y}V_Y\otimes\wedge^{p-v_Y}V_0
\otimes\O_{\bar X}
(\bar B'+\bar L)\\
\parallel&&\parallel\\
\wedge^{\bar p}\bar V_Y\otimes\wedge^{v_Y-\bar p}V_{Y, K}\otimes
\wedge^{p-v_Y}(V_K\cap V_0)\otimes\O_{\bar X}
(\bar B'+\bar L)&\hookrightarrow& 
\wedge^{\bar p}\bar V_Y\otimes\wedge^{v_Y-\bar p}V_{Y, K}\otimes
\wedge^{p-v_Y}V_0\otimes\O_{\bar X}
(\bar B'+\bar L)
\end{matrix}
}

% admits as direct summand the arrow
%\[ C_{\bar p}(\bar X, \bar B', \bar L, \bar V)\otimes\wedge^{p-\bar p}V_K
%\to\wedge^{\bar p}\bar V_Y\otimes\wedge^{v_Y-\bar p}V_{K,Y}\otimes \wedge^{p-v_Y}
%(V_K\cap V_0)\otimes\O_{\bar X}(\bar B'+\bar L).\]
By \cite{el}, Prop. 4.6,  the upper 
left arrow in \eqref{xbar1} induces a surjection on $H^0$ if $d$
is large enough. Now considering
the $E_1$ spectral sequence associated to the resolution \eqref{map}, we see from the 
vanishings of \cite{el} (cf. \eqref{vanishing})
that $H^{q-1}(\I_{Y/X}(B+L))$,
which is 1-dimensional, injects
into $H^0(\O_{\bar X}(\bar B'+\bar L))$. Therefore the image of the
left vertical map in \eqref{map} on $H^{q-1}$ must contain
$\wedge^{v_Y}V_Y\otimes\wedge^{p-v_Y}(V_0\cap V_K)H^{q-1}(\I_{Y/X}(B+L))$.
Consequently,
\eqspl{}{
\dim K_{p,q}(X, B, L)\geq \binom{\dim(V_0\cap V_K)}{p-v_Y}.
}
%\eqref
Choosing $p>v_Y\sim d^{q-1}$, hence $\dim(V_0)\sim d^n$, 
we can arrange things so that $\dim(V_0\cap V_K)\sim d^n$ while
$p-v_Y\sim d$, hence the above lower bound is exponential in $d$.
On the other hand by the very definition of Koszul cohomology,
its dimension is at most exponential in $d$. Hence we conclude:
%its diemension grows at most exponentially in $d$. Therefore we conclude
%and adjusting $V_0, V_K$ appropriately,
%we conclude
\begin{thm}
For $n=\dim(X)\geq 2$, 
$q\geq 2$ and $p$ asymptotically between $d^{q-1}$ and $d^n$, 
the dimension of $K_{p,q}(X, B; L_d)$ grows exponentially with $d$;
in particular,
\[\lim\limits_{d\to\infty}\dim K_{p,q}(X, B, L_d)=\infty.\]
\end{thm}
%In fact, with appropriate choices (e.g. $\dim(V_0\cap V_K)\sim d^n, n>1,  p-v_Y\sim d$), 
%we can arrange that $\dim K_{p,q}(X, B, L_d)$ grows 
%at least like $\binom{d^n}{d}$, i.e. exponentially with $d$. 
This sheds some light on Problem
7.3 in \cite{el}.
\vfill\eject
\bibliographystyle{plain}
\bibliography{../mybib}
\end{document}

%% file: myheader3.tex
 \DeclareMathOperator{\gr}{Gr}

\def\refp #1.{(\ref{#1})}

\def\sbr #1.{^{[#1]}}
\def\sfl #1.{^{\lfloor #1\rfloor}}

\def\?{{\bf{??}}}

\def\H{\mathcal H}

\def\HH{\mathbb H}

\def\P{\mathbb P}

\def\sym{\text{\rm Sym} }

\def\ls{\vskip.25in}

\def\O{\mathcal O}

\def\g{\mathfrak g}

\def\1/2{\frac{1}{2}}

\def\I{\mathcal{ I}}

\def\2{{[2]}}

\def\nl{\newline}

\def\<{\langle}
\def\>{\rangle}

\def\2{{[2]}}

\def\scl #1.{^{\lceil#1\rceil}}
\def\spr #1.{^{(#1)}}
\def\sbc #1.{^{\{#1\}}}

\def\subpr#1.{_{(#1)}}

\def\beq{\begin{equation*}}
\def\eeq{\end{equation*}}

\def\g3{{\Gamma\spr 3.}}

\newcommand{\eqspl}[2]{
\begin{equation}\label{#1}
\begin{split}
#2\end{split}\end{equation}}
\newcommand{\eqsp}[1]{\begin{equation*}
\begin{split}#1\end{split}\end{equation*}}

\newcommand{\exseq}[3]{
0\to #1\to #2\to #3\to 0
}
\newcommand{\beginalphaenum}{
\begin{enumerate}\renewcommand{\labelenumi}{ }
\item \begin{enumerate}
}

\def\eex{\end{rm}\end{example}}
\newtheorem*{lem*}{Lemma}